\documentclass[11pt]{amsart}
\usepackage{amscd,hyperref}
\usepackage[arrow,matrix]{xy}
\usepackage{graphicx}
\usepackage{amsmath, amsfonts}
\usepackage{color}
\input xypic

\numberwithin{equation}{section}
\theoremstyle{plain}

\newtheorem{lemma}{Lemma}[section]

\theoremstyle{definition}
\newtheorem{defn}[lemma]{Definition}
\newtheorem{rem}[lemma]{Remark}

\numberwithin{equation}{section}

\newcommand{\NN}{\mathbb{N}}

\newcommand{\PP}{\mathbb{P}}

\newcommand{\E}{\mathbb{E}}

\newcommand{\X}{\mathbf{X}}

\renewcommand{\i}{\mathbf{i}}
\renewcommand{\k}{\mathbf{k}}
\renewcommand{\j}{\mathbf{j}}
\newcommand{\x}{\mathbf{x}}

\newcommand{\sss}{\mathbf{s}}

\newcommand{\y}{\mathbf{y}}
\newcommand{\z}{\mathbf{z}}

\newcommand{\bel}[1]{\begin{equation}\label{#1}}

\newcommand{\be}{\begin{equation}}
\newcommand{\ba}{\begin{eqnarray}}
\newcommand{\ea}{\end{eqnarray}}

\newcommand{\qe}{\end{equation}}
\newcommand{\al}{\alpha}
\newcommand{\alv}{\boldsymbol{\alpha}}
\newcommand{\bev}{\boldsymbol{\beta}}
\newcommand{\ld}{\lambda}
\newcommand{\de}{\delta}
\newcommand{\Om}{\Omega}

\newcommand{\suml}{\sum\limits}

\newcommand{\ble}{\left\{}
\newcommand{\bri}{\right\}}

\begin{document}

\title[Moment generating functions for Wright-Fisher]{\large The evolution
  of moment generating functions for the Wright-Fisher model of
  population genetics}
\author{Tat Dat Tran, Julian Hofrichter, J\"{u}rgen Jost}
\date{\today}

\medskip

\abstract We derive and apply a partial differential equation for the
moment generating function of the Wright-Fisher model of population
genetics.  

\endabstract

\maketitle
\tableofcontents

{\it  MSC2000: 60J27, 60J60}

{\it Key words: Wright-Fisher model, random genetic drift, moment generating function, master  equation}

\section{Introduction}\label{sec:intro}
The Wright-Fisher model, that is, the random genetic drift model
developed  by Fisher  \cite{fisher} and  by Wright  \cite{wright1} and
mathematically solved by Kimura \cite{kimura1,kimura2} is
the basic stochastic model in population genetics
(see for instance \cite{ewens}). The discrete model  is concerned with the
evolution of the probabilities between non-overlapping generations in
a population of fixed size of two or more alleles obtained from
random sampling in the parental generation. This basic model thus
describes random genetric drift, and  additional biological mechanisms
like mutation, selection, or a spatial population structure can then
be superposed. As such, the model works with a finite population in
discrete time, but the mathematical analysis of Kimura and others
turned to its diffusion approximation suggested by Kolmogorov. This
diffusion approximation works with an infinite population in
continuous time. It consists of two partial differential equations of
parabolic type for the probability density function for the various
alleles, the so-called forward and backward Kolmogorov equations.  In  \cite{thj1,thj2}, we have presented a general
solution scheme for the associated diffusion process that keeps track
of the population across possible allele losses. With our scheme, all
basic quantities of interest, like expected times of allele losses,
can be readily derived. In other words, we derive a global solution,
in contrast to the local ones of Kimura and others. A crucial
ingredient in our scheme are the equations for the moments of the
probability distribution. 

A somewhat simpler model than the Wright-Fisher model that however leads to the same diffusion
approximation is the Moran model, a simple birth-death process in
continuous time, see \cite{ewens}. In the Wright-Fisher model, when
creating the next generation, for each new  member of the population
a parent in the previous generation is randomly chosen. This has the consequence that one and the same individual in that parent
generation could produce several offspring. In the Moran model, in
contrast, a
randomly chosen individual gives birth to a clone, and then another
random 
individual in the population is killed to keep the population size
constant. Thus, here, in each step, only one offspring is produced. 

As mentioned, the Kolmogorov diffusion equations are concerned with
asymptotic quantities, and in particular, do not account for small
population size effects. Therefore, Houchmandzadeh and
Vallade \cite{HV2010} have proposed to use the master equation for the
probability distribution to directly derive a partial differential
equation for the probability generating function of the process. This
approach can produce exact formulae even for finite populations. In
\cite{HV2010}, this has been carried out for the Moran model with two alleles, not only for
the basic model of random genetic drift, but also including
the case of selection.

In the present paper, we derive a partial differential equation for
the exponential moment generating function of the Wright-Fisher model
with arbitrarily many alleles,
utilizing our scheme of moment equations mentioned above. From this
scheme, we can then also easily rederive formulae for quantities of
interest like fixation probabilities. In order to facilitate the
understanding, we shall always first treat the simplest case of two
alleles and then present the case of arbitrarily many alleles.

\section{Master equation}\label{sec:master equation}
In this section we shall use the master equations that express the
evolution of the probability distribution for the alleles in the
population in terms of the transition probabilities, in order to
derive differential equations for the moments of the process. 

%%%%%%%%%%%%%%%%%%%%%%%%%%%%%%%%%%%%%%
%%%%%%%%   Subsection 2.1
\subsection{2 alleles}
Consider a continuous time stochastic process $\{X_t\}_{t\ge 0}$ with values in 
$$
S_{1}^{2N}=\Bigg\{0, \frac{1}{2N},\ldots,1\Bigg\},
$$
with transition rates $B(k,j)$ from state $k$ to $j$ specified below. 

The master equation for the  probability function 
$$P(t,i,j)=\PP\Bigg(X_t=\frac{j}{2N}\Bigg|X_0=\frac{i}{2N}\Bigg)$$ will be
\bel{eq:master}
\frac{\partial P(t,i,j)}{\partial t} = \suml_{k=0}^{2N}P(t,i,k) B(k,j), \quad \forall i,j=\overline{0,2N},
\qe
with the initial values $P(0,i,j)=\delta_{ij}$. In  matrix form, this
reads as 
\bel{eq:matrix}
\begin{cases}
\frac{\partial P(t)}{\partial t} =& P(t) B,\quad \forall t\ge 0\\
P(0) =& I
\end{cases}
\qe
This is a linear problem (\ref{eq:matrix}), with the unique solution
$P(t)=e^{Bt}$. $B$ and $P(t)$ then commute, and therefore we also have
\bel{eq:master2}
\frac{\partial P(t,i,j)}{\partial t} = \suml_{k=0}^{2N}P(t,k,j) B(i,k), \quad \forall j=\overline{0,2N}.
\qe
For $B$, we consider two cases
\bel{formula-a}
\sum_{k=0}^{2N} \Big(\frac{k}{2N}\Big)^n B(i,k)=\frac{n(n-1)}{2}\Bigg(\Big(\frac{i}{2N}\Big)^{n-1}-\Big(\frac{i}{2N}\Big)^{n}\Bigg),\quad \forall n\ge0; \: i=\overline{0,2N}
\qe
and
\bel{formula-b}
B(i,k)=\binom{2N}{k}\Big(\frac{i}{2N}\Big)^k \Big(1-\frac{i}{2N}\Big)^{2N-k}-\de_{i,k},\quad  \: i,k=\overline{0,2N}
\qe

\begin{rem}
\begin{itemize}
\item Case \eqref{formula-b} corresponds to the Wright-Fisher model,
  more precisely a continuous time Wright-Fisher model with discrete
  states, whereas
the implicit scheme \eqref{formula-a} will get rid off certain error
terms for finite population size. 
\item In  \eqref{formula-a}, which has been directly constructed from the moment equation, the coefficients $B(i,k)$ could become
  negative (see for example \cite{Gar2004} for such a generalization
  of the master equation concept).
\item Since in general the transition rates $B(i,j)$ are nonzero for
  any pair $i,j$,  in our master equations, starting from state $\i$
  we can directly access any other state with positive
  probability. Therefore, for  the probability generating function, we
  would get an $2N-$order partial differential equation which may be
  hard to solve. In contrast, for  the moment generating function, we
  shall get a second order partial differential equation which can be
  solved by a simple expansion.
\end{itemize}
\end{rem}

We shall prove that in the limit of the population size ($2N\to
\infty$), these master equations will produce the classical
Wright-Fisher diffusion equations. Moreover we shall prove that these
master equations satisfy the moment formulae, exactly ofr
\eqref{formula-a} and approximately for \eqref{formula-b}. Therefore we can apply the moment generating function technique to calculate the  conditional probability function.

In fact, when $2N$ is sufficient large, we set
$$
x=\frac{i}{2N},\quad y=\frac{j}{2N},\quad z=\frac{k}{2N},\quad dx=dy=dz=\frac{1}{2N}
$$
and
$$
p(t,x,y)dy=P(t,i,j), \quad b(x,y)=B(i,j).
$$
Then, we obtain from (\ref{eq:master2})
$$
\frac{\partial p(t,x,y)}{\partial t} = \suml_{z} p(t,z,y) b(x,z), \quad \forall y\in S_1^{2N},
$$
Now we expand the function $p(t,z,y)$ in $z$ at $x$ and obtain
$$
p(t,z,y)=\suml_{n\ge 0}\frac{1}{n!}\frac{\partial^n p(t,x,y)}{\partial x^n} (z-x)^n.
$$
\begin{enumerate}
\item
In the case of  (\ref{formula-a}), we obtain 
$$
\sum_{z} z^n b(x,z)=\frac{n(n-1)}{2}\Bigg(x^{n-1}-x^{n}\Bigg),\quad \forall n\ge0; \: x\in S_1^{2N}.
$$
It follows by induction that
\bel{eq:b}
\begin{cases}
\suml_{z} (z-x)^n b(x,z) =& 0, \quad n\ne 2\\
\suml_{z} (z-x)^2 b(x,z) =& x(1-x).
\end{cases}
\qe
Then we obtain
\be
\begin{split}
\frac{\partial p(t,x,y)}{\partial t} =& \suml_{z} p(t,z,y) b(x,z)\\
=&\suml_{z} \Bigg\{\suml_{n\ge 0}\frac{1}{n!}\frac{\partial^n p(t,x,y)}{\partial x^n} (z-x)^n\Bigg\} b(x,z)\\
=& \frac{x(1-x)}{2} \frac{\partial^2 p(t,x,y)}{\partial x^2}
\end{split}
\qe
which is exactly the classical Wright-Fisher diffusion equation.
\item
In the case of  (\ref{formula-b}), we obtain 
$$
\sum_{z} z^n b(x,z)=\frac{n(n-1)}{4N}\Bigg(x^{n-1}-x^{n}\Bigg)+O\Big(\frac{1}{N^2}\Big),\quad \forall n\ge0; \: x\in S_1^{2N}.
$$
It follows by induction that
\bel{eq:b}
\begin{cases}
\suml_{z} (z-x)^n b(x,z) =& O\Big(\frac{1}{N^2}\Big), \quad n\ne 2\\
\suml_{z} (z-x)^2 b(x,z) =& \frac{x(1-x)}{2N}.
\end{cases}
\qe
Then we obtain
\be
\begin{split}
\frac{\partial p(t,x,y)}{\partial t} =& \suml_{z} p(t,z,y) b(x,z)\\
=&\suml_{z} \Bigg\{\suml_{n\ge 0}\frac{1}{n!}\frac{\partial^n p(t,x,y)}{\partial x^n} (z-x)^n\Bigg\} b(x,z)\\
=& \frac{x(1-x)}{4N} \frac{\partial^2 p(t,x,y)}{\partial x^2}+O\Big(\frac{1}{N^2}\Big)
\end{split}
\qe
which is an approximation of the  Wright-Fisher diffusion equation.
\end{enumerate}

Now we shall use the master equation (\ref{eq:master}) to derive the moment equation.

In fact, the $n-th$ moment of this conditional probability function is 
$$
m_n(t):=\suml_{j=0}^{2N} \Bigg(\frac{j}{2N}\Bigg)^n P(t,i,j).
$$
Therefore we have
\bel{eq:moment}
\begin{split}
\dot{m}_n(t)=&\suml_{y} y^n \frac{\partial p(t,x,y)}{\partial t} \frac{1}{2N}\\
=& \suml_{y} y^n\suml_{z} p(t,x,z) b(z,y) \frac{1}{2N} \quad\text{by the Master equation (\ref{eq:master})}\\
=& \suml_{z} \Bigg(\suml_{y} y^n b(z,y)\Bigg) p(t,x,z) \frac{1}{2N}\\
=& \suml_{z} \frac{n(n-1)}{2}(z^{n-1}-z^n)  p(t,x,z) \frac{1}{2N}
\quad\text{by the formula for $b$ from (\ref{formula-a})}\\
&\Bigg( =\suml_{z} \frac{n(n-1)}{4N}(z^{n-1}-z^n)  p(t,x,z)
\frac{1}{2N}+O\Big(\frac{1}{N^2}\Big) \quad\text{by the formula for
  $b$ from (\ref{formula-b})}\Bigg)\\
=&-\frac{n(n-1)}{2} m_n(t)+\frac{n(n-1)}{2} m_{n-1}(t) \quad\text{in the case of (\ref{formula-a})}\\
&\Bigg(=-\frac{n(n-1)}{4N} m_n(t)+\frac{n(n-1)}{4N} m_{n-1}(t)+O\Big(\frac{1}{N^2}\Big) \quad\text{in the case of (\ref{formula-b})}\Bigg).
\end{split}
\qe

%%%%%%%%%%%%%%%%%%%%%%%%%%%%%%%%%%%%%%
%%%%%%%%   Subsection 2.2
\subsection{$K+1$ alleles}

Consider a continuous time stochastic process $\{\X_t\}_{t\ge 0}$ with values in 
$$
S_{K}^{2N}:=\Bigg\{\frac{\i}{2N}=\Big(\frac{i_1}{2N},\cdots,\frac{i_K}{2N}\Big): i_u \in \NN_0 \:\text{ for all }  u=\overline{1,K},\: \text{and } \suml_{u=1}^K i_u \le 2N\Bigg\}.
$$
To simplify the notation, we also put 
$$
\Om_{K}^{2N}:=2N S^{2N}_K=\Bigg\{\i=\Big(i_1,\cdots,i_K\Big): i_u \in \NN_0 \:\text{ for all }  u=\overline{1,K},\: \text{and } \suml_{u=1}^K i_u \le 2N\Bigg\}.
$$ 

The master equation for the conditional probability function 
$$P(t,\i,\j)=\PP\Bigg(\X_t=\frac{\j}{2N}\Bigg|\X_0=\frac{\i}{2N}\Bigg)$$
then is 
\bel{eq:Master}
\frac{\partial P(t,\i,\j)}{\partial t} = \suml_{\k\in \Om^{2N}_K} P(t,\i,\k) B(\k,\j), \quad \forall \i,\j\in \Om^{2N}_K,
\qe
with  initial values $P(0,\i,\j)=\delta_{\i\j}$. In  matrix form,  we have
\bel{eq:Matrix}
\begin{cases}
\frac{\partial P(t)}{\partial t} =& P(t) B,\quad \forall t\ge 0\\
P(0) =& I
\end{cases}
\qe
As for 2 alleles, the problem (\ref{eq:Matrix}) has a unique solution
$P(t)=e^{Bt}$, and $B$ and $P(t)$ commute, and therefore we also have
\bel{eq:Master2}
\frac{\partial P(t,\i,\j)}{\partial t} = \suml_{\k\in \Om^{2N}_K}P(t,\k,\j) B(\i,\k), \quad \forall \i,\j\in \Om^{2N}_K,
\qe
where $B$ is  defined by the formulae
\bel{Formula}
\begin{split}
\sum_{k\in \Om^{2N}_K} \Big(\frac{\k}{2N}\Big)^{\alv} B(\i,\k)=-\frac{|\alv|(|\alv|-1)}{2}\Big(\frac{\i}{2N}\Big)^{\alv}-\suml_{u=1}^K \frac{\al_u (\al_u-1)}{2}\Big(\frac{\i}{2N}\Big)^{\alv-e_u},\\
\quad \forall \alv\in \NN_0^K; \: \i\in \Om^{2N}_K.
\end{split}
\qe

We shall prove that in the limit of the population  size ($2N\to
\infty$), this master equation will yield the  Wright-Fisher diffusion
equation. Again, this master equation will satisfy the moment formulae
and  we can apply the moment generating function technique to calculate the  conditionalal probability function.

In fact, when $2N$ is sufficient large, we set
$$
\x=\frac{\i}{2N},\quad \y=\frac{\j}{2N},\quad \z=\frac{\k}{2N},\quad dx_u=dy_u=dz_u=\frac{1}{2N}
$$
and
$$
p(t,\x,\y)\frac{1}{2N}=P(t,\i,\j),\quad \quad b(\x,\y)=B(\i,\j).
$$
Then, we obtain from  (\ref{eq:Master2})
$$
\frac{\partial p(t,\x,\y)}{\partial t} = \suml_{\z} p(t,\z,\y) b(\x,\z), \quad \forall \y\in S_K^{2N},
$$
Now we expand the function $p(t,\z,\y)$ in $\z$ at $\x$ to obtain
$$
p(t,\z,\y)=\suml_{\alv}\frac{1}{\alv!}\frac{\partial^{\alv} p(t,\x,\y)}{\partial \x^{\alv}} (\z-\x)^{\alv}.
$$
We obtain from the formulae (\ref{Formula}) that
$$
\sum_{\z} \z^{\alv} b(\x,\z)=-\frac{|\alv|(|\alv|-1)}{2}\x^{\alv} + \sum_{u=1}^K \frac{\al_u(\al_u-1)}{2} \x^{\alv-e_u},\quad \forall \alv\in\NN_0^K; \: \x\in S_K^{2N}.
$$
It follows by induction that
\bel{eq:B}
\begin{cases}
\suml_{\z} (\z-\x)^{\alv} b(\x,\z) =& 0, \quad |\alv|\ne 2\\
\suml_{\z} (\z-\x)^{e_u+e_v} b(\x,\z) =& x_u(\de_{uv}-x_v), \: u,v=\overline{1,K}.
\end{cases}
\qe
Then we obtain
\be
\begin{split}
\frac{\partial p(t,\x,\y)}{\partial t} =& \suml_{\z} p(t,\z,\y) b(\x,\z)\\
=&\suml_{\z} \Bigg\{\suml_{\alv\in\NN_0^K}\frac{1}{\alv!}\frac{\partial^{\alv} p(t,\x,\y)}{\partial \x^{\alv}} (\z-\x)^{\alv}\Bigg\} b(\x,\z)\\
=& \suml_{u,v=1}^K\frac{x_u(\de_{uv}-x_v)}{2} \frac{\partial^2 p(t,\x,\y)}{\partial x_u \partial x_v}
\end{split}
\qe
which is exactly the  Wright-Fisher diffusion equation for $K+1$ alleles.

Now we shall prove that the master equation (\ref{eq:master}) yields the moment equations.

In fact, the $\alv-th$ moment of this conditional probability function is 
$$
m_{\alv}(t):=\suml_{\j\in \Om^{2N}_K} \Bigg(\frac{\j}{2N}\Bigg)^{\alv} P(t,\i,\j).
$$
Therefore we have
\bel{eq:moment2}
\begin{split}
\dot{m}_{\alv}(t)=&\suml_{\y} \y^{\alv} \frac{\partial p(t,\x,\y)}{\partial t} \frac{1}{2N}\\
=& \suml_{\y} \y^{\alv}\suml_{\z} p(t,\x,\z) b(\z,\y) \frac{1}{2N} \quad\text{due to the Master equation (\ref{eq:Master})}\\
=& \suml_{\z} \Bigg(\suml_{\y} \y^{\alv} b(\z,\y)\Bigg) p(t,\x,\z) \frac{1}{2N}\\
=& \suml_{\z} \Bigg(-\frac{|\alv|(|\alv|-1)}{2}\z^{\alv}+\suml_{u=1}^K \frac{\al_u(\al_u-1)}{2} \z^{\alv-e_u})\Bigg)  p(t,\x,\z) \frac{1}{2N}\\
& \quad\qquad \text{by the formulae for $b$ from (\ref{Formula})}\\
=&-\frac{|\alv|(|\alv|-1)}{2} m_{\alv}(t)+\suml_{u=1}^K \frac{\al_u(\al_u-1)}{2} m_{\alv-e_u}(t).
\end{split}
\qe

%%%%%%%%%%%%%%%%%%%%%%%%%%%%%%%%%%%%%%
%%%%%%%%   Section 3
\section{Moment generating functions}\label{sec:MGF}

\begin{defn} 
\begin{enumerate}
\item Let $X$ be a random variable with discrete values with probability distribution function $p(x)=\PP\big[X=x\big]$. The (exponential) moment generating function of the random  variable $X$ is
$$
H(s)  := \E\Big[ e^{sX}\Big] = \suml_{x} e^{xs} p(x)                            
$$
(defined for those values of $s \in {\mathbb R}$ for which the sum converges).
\item Let $\X=(X^1,\cdots,X^K)$ be a tuple of random variables with the joint probability distribution function $p(x^1,\ldots,x^K)=\PP\Big[X^1=x^1,\ldots,X^K=x^K\Big]$. The (exponential) moment generating function of $\X$ then is
$$
H(s_1,\ldots,s_K)  := \E\Big[ e^{\suml_{i=1}^K s_i X^i}\Big] = \suml_{x^1,\ldots, x^K} e^{\suml_{i} s_i x^i} p(x^1,\ldots,x^K)=\suml_{\x} e^{\sss \cdot \x} p(\x).                            
$$
(defined for those values of $\sss \in {\mathbb R^K}$ for which the sum converges).
\end{enumerate}
\end{defn}

Here, the moments of $\X$ can directly be computed from the
derivatives of $H(\sss)$ at $\sss = 0$,  
$$
\E\Big[\X^{\alv}\Big]=\frac{\partial^{\alv}  H(\sss)}{\partial \sss^{\alv}}{\Big|_{\sss=0}}.
$$
We shall now derive the (second order) partial differential equation for the (exponential) moment generating functions of our Markov process $\X_t$.
%%%%%%%%%%%%%%%%%%%%%%%%%%%%%%%%%%%%%%
%%%%%%%%   Subsection 3.1
\subsection{2 alleles}

The exponential moment generating function is
\be
\begin{split}
H(t;s) = & \E[e^{sX_t}]\\
=&\suml_{n\ge 0} \frac{s^n}{n!}\E[(X_t)^n]\\
=&\suml_{n\ge 0} \frac{s^n}{n!}m_n(t),
\end{split}
\qe
where $m_n(t)$ is the $n-th$ moment of $X_t$ around $0$.

From the equation \eqref{eq:moment} for the moments
$$ 
\dot{m}_n(t)=-\frac{n(n-1)}{2}m_n(t)+\frac{n(n-1)}{2}m_{n-1}(t)
$$
we obtain
\bel{eq:1}
\begin{split}
\frac{\partial H(t;s)}{\partial t} =&\suml_{n\ge 0}\dot{m}_n(t) \frac{s^n}{n!}\\
=&\suml_{n\ge 0}\Bigg[-\frac{n(n-1)}{2}m_n(t)+\frac{n(n-1)}{2}m_{n-1}(t)\Bigg] \frac{s^n}{n!}\\
=&\suml_{n\ge 2}-\frac{1}{2}m_n(t) \frac{s^n}{(n-2)!}+\suml_{n\ge 2}\frac{1}{2}m_{n-1}(t) \frac{s^n}{(n-2)!}\\
=&-\frac{s^2}{2}\suml_{n\ge 0}m_{n+2}(t)\frac{s^n}{n!} +\frac{s^2}{2}\suml_{n\ge 0}m_{n+1}(t)\frac{s^n}{n!}\\
=&-\frac{s^2}{2}\frac{\partial ^2}{\partial s^2} H(t;s)+\frac{s^2}{2}\frac{\partial}{\partial s} H(t;s).
\end{split}
\qe

We now consider solutions of such equations. First, we solve the
equation (\ref{eq:1}) by  separation of variables. With
$H(t;s)=T(t)S(s)$, the equation becomes
$$
\frac{T'(t)}{T(t)}=\frac{-s^2S''(s)+s^2 S'(s)}{2S}=-\ld.
$$ 
It follows that $T(t)=Ce^{-\ld t}$ and $S(s)$ satisfies the ODE
\bel{eq:2}
-x^2 y_{xx} +x^2 y_x = -2\ld y.
\qe
By putting $y(x)=\suml_{n\ge 0}a_n x^n$ and equating coefficients in
the  ODE (\ref{eq:2}) we obtain:
\begin{enumerate}
\item If 
$$
\ld \not\in \Lambda:=\Big\{\mu_n=\frac{n(n-1)}{2}, n\in \NN\Big\}
$$
then the ODE (\ref{eq:2}) has a unique solution $y(x)= 0$;
\item If $\ld=\mu_0$ then $y_0(x)=a^{(0)}_0:=1$;
\item If $\ld = \mu_k$ for some $k\ge 1$ then the solution is of the form
\be
y_k(x)=\suml_{n\ge 0}a^{(k)}_n x^n
\qe
where
\bel{eq:3}
a^{(k)}_n=
\begin{cases}
&0, \quad \text{if } n<k\\
&1, \quad \text{if } n=k\\
&\frac{n-1}{2(\mu_n-\mu_k)} \cdots  \frac{k}{2(\mu_{k+1}-\mu_k)},\quad \text{if } n\ge k+1.
\end{cases}
\qe
\end{enumerate}
Therefore the solution of  (\ref{eq:1}) is 
\be
\begin{split}
H(t;s)=&\suml_{k\ge 0} c_k y_k(s) e^{-\mu_k t}\\
=&\suml_{k\ge 0} c_k \Bigg( \suml_{n\ge k} a^{(k)}_n s^n\Bigg) e^{-\mu_k t}\\
=&\suml_{n\ge 0} \Bigg(n!\suml_{k=0}^n c_k a^{(k)}_n e^{-\mu_k t}\Bigg) \frac{s^n}{n!}.
\end{split}
\qe
This yields the  moment formula
$$
m_n(t)=n!\suml_{k=0}^n c_k a^{(k)}_n e^{-\mu_k t}=\suml_{k=0}^n c_k A^{(k)}_n e^{-\mu_k t}.
$$
The coefficients $c_k$ can be  calculated from  the initial condition
$$
\Big(\frac{i}{2N}\Big)^n:=p^n=m_n(0)=\suml_{k=0}^n c_k A^{(k)}_n,\quad \forall n\ge 0.
$$
In fact, by representing these equalities in  matrix form
\be
\begin{bmatrix}
1 & 0 &\cdots & 0&0\\
A^{(0)}_1 & 1 & \cdots &0&0\\
\vdots & \vdots & \ddots & \vdots\\
A^{(0)}_{n-1} & A^{(1)}_{n-1}&\cdots &(n-1)! &0\\
A^{(0)}_{n} & A^{(1)}_{n}&\cdots &A^{(n-1)}_{n} &n!
\end{bmatrix}
\begin{bmatrix}
c_0\\
c_1\\
\vdots \\
c_{n-1}\\
c_n
\end{bmatrix}
=
\begin{bmatrix}
1\\
p\\
\vdots\\
p^{n-1}\\
p^{n}
\end{bmatrix}, 
\qe
it follows that
\be
\begin{bmatrix}
c_0\\
c_1\\
\vdots \\
c_{n-1}\\
c_n
\end{bmatrix}
=\begin{bmatrix}
1 & 0 &\cdots & 0&0\\
A^{(0)}_1 & 1 & \cdots &0&0\\
\vdots & \vdots & \ddots & \vdots\\
A^{(0)}_{n-1} & A^{(1)}_{n-1}&\cdots &(n-1)! &0\\
A^{(0)}_{n} & A^{(1)}_{n}&\cdots &A^{(n-1)}_{n} &n!
\end{bmatrix}^{-1}
\begin{bmatrix}
1\\
p\\
\vdots\\
p^{n-1}\\
p^{n}
\end{bmatrix}. 
\qe
\begin{rem}
We can easily check some instances:
Because of $c_0=1; c_1=p; c_2=\frac{p^2-p}{2}; c_3=\frac{p^3-3/2 p^2+1/2 p}{6}$ then $m_0(t)=1; m_1(t)=p; m_2(t)=p+(p^2-p)e^{-t}; m_3(t)=p+3/2(p^2-p)e^{-t}+(p^3-3/2 p^2+1/2 p)e^{-3t}$.
\end{rem}

This also yields  the fixation probability at time $t$ is (see also  \cite{thj1})
\be
\begin{split}
P(t,i,2N)=&\lim_{n\to \infty}  m_n(t)\\
=&\lim_{n\to \infty} \suml_{k=0}^n c_k A^{(k)}_n e^{-\mu_k t}\\
=&p+\lim_{n\to \infty} \suml_{k=2}^n c_k A^{(k)}_n e^{-\mu_k t},
\end{split}
\qe
and the eventual fixation probability
\be
\begin{split}
P(\infty,i,2N)=&\lim_{t\to \infty}\Bigg(p+\lim_{n\to \infty} \suml_{k=2}^n c_k A^{(k)}_n e^{-\mu_k t}\Bigg)\\
=&p+\lim_{n\to \infty} \lim_{t\to \infty} \suml_{k=2}^n c_k A^{(k)}_n e^{-\mu_k t}\\
=&p.
\end{split}
\qe

Similarly, by calculating for the other allele $(Y_t=1-X_t)$, we obtain the extinction probability at time $t$ is (see also  \cite{thj1})
\be
\begin{split}
P(t,i,0)=&\lim_{n\to \infty}  m'_n(t)\\
=&1-p+\lim_{n\to \infty} \suml_{k=0}^n c'_k A^{(k)}_n e^{-\mu_k t}\\
=&1-p+\lim_{n\to \infty} \suml_{k=2}^n c'_k A^{(k)}_n e^{-\mu_k t},
\end{split}
\qe
and the eventual extinction probability
\be
\begin{split}
P(\infty,i,0)=&\lim_{t\to \infty}\Bigg(1-p+\lim_{n\to \infty} \suml_{k=2}^n c_k A^{(k)}_n e^{-\mu_k t}\Bigg)\\
=&1-p.
\end{split}
\qe

The moments of the sojourn and absorption times  were derived by Nagylaki \cite{nagylaki} for two alleles, and by Lessard and Lahaie \cite{lessard} in the multi-allele case. We denote by $T^{k+1}_{n+1}(p)=\inf \ble{t>0: X_t\in \overline{V}_k}|X_0=p \bri$ the first time when the  population has (at most) $k+1$ alleles. $T^{k+1}_{n+1}(p)$ is a continuous random variable valued in $[0,\infty)$ and we denote by $\phi(t,p)$ its probability density function. It is easy to see that $\overline{V}_k$ is invariant under the process $(X_t)_{t \ge 0}$, i.e. if $X_s \in \overline{V}_k$ then $X_t\in \overline{V}_k$ for all $t\ge s$ (once an allele is lost from the population, it can never again be recovered). We have the equality
$$
\PP(T^{1}_{2}(p)\le t)=P(t,i,0)+P(t,i,2N).
$$
Therefore the expectation of the absorption time is
\begin{equation*}
\begin{split}
\E(T^{1}_{2}(p))=&\int_0^\infty t \frac{\partial }{\partial t}\Big(P(t,i,0)+P(t,i,2N)\Big)dt\\
=&-\lim_{n\to \infty} \sum_{k=2}^n (c_k+c'_k) A_n^{(k)}\frac{1}{\mu_k}.
\end{split}
\end{equation*}

Moreover we have
\be
\suml_{j=0}^{2N} \Bigg(\frac{j}{2N}\Bigg)^n P(\infty,i,j)=\lim_{t\to \infty}  m_n(t)=
\begin{cases}
p,& \text{ for } n\ge 1\\
1,& \text{ for } n= 0
\end{cases}
\qe
Therefore we obtain the eventually probability
$$
P(\infty,i,j)=p\de_{2N,j}+(1-p)\de_{0,j}.
$$

The probability of heterogeneity is (also see \cite{thj1})
\be
\begin{split}
H_t:=&2\suml_{j=0}^{2N} \frac{j}{2N} \Big(1-\frac{j}{2N} \Big) P(t,i,j)\\
=&2(m_1(t)-m_2(t))\\
=&2\Bigg(p-\Big(p+c_2 A_2^{(2)}e^{-t}\Big)\Bigg)\\
=&2p(1-p)e^{-t}.
\end{split}
\qe
%%%%%%%%%%%%%%%%%%%%%%%%%%%%%%%%%%%%%%
%%%%%%%%   Subsection 3.2 
\subsection{$K+1$ alleles}

We can apply the same scheme for any $K$. The exponential generating
function now is 
$$
H(t;s_1,\ldots,s_K)=\suml_{\alv}m_{\alv}(t) \frac{\sss^{\alv}}{\alv!},
$$
where $m_{\alv}(t)$ is the $\alv-th$ moment of $\X_t$ around $0$.

From \eqref{eq:moment2}, i.e., 
\bel{cond:moment}
\dot{m}_{\alv}(t)=-\frac{|\alv|(|\alv|-1)}{2}
m_{\alv}(t)+\suml_{i=1}^K \frac{\al_i (\al_i-1)}{2}m_{\alv-e_i}(t), 
\qe

we obtain
\bel{eq:4}
\begin{split}
\frac{\partial H(t;\sss)}{\partial t} =&\suml_{\alv}\dot{m}_{\alv}(t) \frac{\sss^{\alv}}{\alv!}\\
=&\suml_{\alv}\Bigg[-\frac{|\alv|(|\alv|-1)}{2} m_{\alv}(t)+\suml_{i=1}^K \frac{\al_i (\al_i-1)}{2}m_{\alv-e_i}(t)\Bigg] \frac{\sss^{\alv}}{\alv!}\\
=&\suml_{\alv}\Bigg(-\frac{\sum_{i\ne j} \al_i\al_j }{2}-\frac{\sum_{i} \al_i(\al_i-1) }{2}\Bigg)m_{\alv}(t) \frac{\sss^{\alv}}{\alv!}\\
&+\suml_{\alv}\suml_{i=1}^K \frac{\al_i (\al_i-1)}{2}m_{\alv-e_i}(t)\frac{\sss^{\alv}}{\alv!}\\
=&-\frac{1}{2}\suml_{i\ne j}  s_i s_j \frac{\partial ^2 H(t,\sss)}{\partial s_i \partial s_j} -\frac{1}{2} \suml_{i}s_i^2 \frac{\partial^2 H(t,\sss)}{\partial s_i^2}+\suml_{i=1}^K \frac{1}{2} \suml_{i}s_i^2 \frac{\partial H(t,\sss)}{\partial s_i} \\
=&-\frac{1}{2}\suml_{i,j=1}^K s_i s_j \frac{\partial ^2}{\partial s_i \partial s_j}H(t;\sss)+\suml_{i=1}^K \frac{s_i^2}{2} \frac{\partial }{\partial s_i}H(t;\sss)
\end{split}
\qe

Separating variables  as above,  $T(t)=Ce^{-\ld t}$ and $S(\sss)$ satisfies the PDE
\bel{eq:5}
-\frac{1}{2}\suml_{i,j=1}^K s_i s_j \frac{\partial ^2}{\partial s_i \partial s_j}y(\sss)+\suml_{i=1}^K \frac{s_i^2}{2} \frac{\partial }{\partial s_i} y(\sss)=-\ld y(\sss).
\qe

By putting $y(\sss)=\suml_{\alv} a_{\alv} {\sss}^{\alv}$ and equating
coefficients in the  PDE (\ref{eq:5}) we obtain:
\begin{enumerate}
\item If 
$$
\ld \not\in \Lambda:=\Big\{\mu_n=\frac{n(n-1)}{2}, n\in \NN\Big\}
$$
then the PDE (\ref{eq:5}) has a unique solution $y(\sss)= 0$;
\item If $\ld=\mu_0$ then $y_0(x)=a^{(0)}_0:=1$;
\item If $\ld = \mu_k$ for some $k\ge 1$ then there are $\binom{k}{2}$ independent solutions of the form
\be
y_{k,\alv}(\sss)=\suml_{\bev}a^{(k)}_{\alv,\bev} \sss^{\bev}, \quad \forall |\alv|=k.
\qe
where
\be
a^{(k)}_{\alv,\bev}=
\begin{cases}
&0, \quad \text{if }  |\bev|< k\\
&\de^{\alv}_{\bev}, \quad \text{if } |\bev|=k\\
&\text{inductively  defined by  (\ref{eq:if}) below},\quad \text{if } |\bev|\ge k+1
\end{cases}
\qe
\bel{eq:if}
a^{(k)}_{\alv,\bev}=\frac{\suml_{i=1}^K (\beta_i-1)a^{(k)}_{\alv,\bev-e_i}}{|\bev|(|\bev|-1)-k(k-1)}.
\qe

Therefore the solution of equation (\ref{eq:5}) is 
\be
\begin{split}
H(t;\sss)=&\suml_{k\ge 0} \suml_{|\alv|=k} c_{k,\alv}y_{k,\alv}(\sss) e^{-\mu_k t}\\
=&\suml_{k\ge 0} \suml_{|\alv|=k} c_{k,\alv} \Bigg( \suml_{\bev} a^{(k)}_{\alv,\bev} \sss^{\bev} \Bigg)e^{-\mu_k t}\\
=&\suml_{\bev} \Bigg(\suml_{k=0}^{|\bev|} \suml_{|\alv|=k} c_{k,\alv}  a^{(k)}_{\alv,\bev} e^{-\mu_k t}\Bigg) \sss^{\bev}\\
=&\suml_{\bev} \bev!\Bigg(\suml_{|\alv|\le |\bev|} c_{|\alv|,\alv}  a^{(|\alv|)}_{\alv,\bev} e^{-\mu_k t}\Bigg) \frac{\sss^{\bev}}{\bev!}
\end{split}
\qe
This yields the  moment formula 
$$
m_{\bev}(t)=\bev!\suml_{|\alv|\le |\bev|} c_{|\alv|,\alv}  a^{(|\alv|)}_{\alv,\bev} e^{-\mu_k t}
$$
where the coefficients $c_{|\alv|,\alv}$ can be computed from  the initial condition
$$
p^{\bev}=m_{\bev}(0)=\bev!\suml_{|\alv|\le |\bev|} c_{|\alv|,\alv}  a^{(|\alv|)}_{\alv,\bev} ,\quad \forall \bev.
$$

\end{enumerate}

\section*{Acknowledgement}  The research leading to these results has received funding from the European Research Council under the European Union Seventh Framework Programme (FP7/2007-2013)/ERC grant agreement no. 267087.

\bigskip
\author{Tat Dat Tran},
\address{Max Planck Institute for Mathematics in the Sciences, Inselstrasse 22, 04103 Leipzig, Germany, }
\email{trandat@mis.mpg.de}\\

\author{Julian Hofrichter},
\address{Max Planck Institute for Mathematics in the Sciences, Inselstrasse 22, 04103 Leipzig, Germany, }
\email{julian.hofrichter@mis.mpg.de}\\

\author{J\"{u}rgen Jost},
\address{Max Planck Institute for Mathematics in the Sciences, Inselstrasse 22, 04103 Leipzig, Germany\\

Santa Fe Institute for the Sciences of Complexity, Santa Fe, NM 87501, USA, }\email{jost@mis.mpg.de}\\

\end{document}